\documentclass[journal]{IEEEtran}   % The CDC formate.

\usepackage{graphics} % for pdf, bitmapped graphics files
\usepackage{epsfig} % for postscript graphics files
\usepackage{amsmath} % assumes amsmath package installed
\usepackage{amssymb}  % assumes amsmath package installed
\usepackage{graphicx}
\usepackage{cite}
\usepackage{amsmath}
\usepackage{stfloats}
\usepackage{balance}
\usepackage{color}
\usepackage{amsfonts}
\usepackage{amssymb}
\usepackage[ruled,linesnumbered]{algorithm2e}
\allowdisplaybreaks
\begin{document}  
	\title{Data-Driven Predictive Control for Continuous-Time Industrial Processes with Completely Unknown Dynamics}
	\author{Yuanqiang Zhou, Dewei Li, Yugeng Xi
		\thanks{Y.~Zhou, D.~Li, Y. G. ~Xi  are with the Department of Automation, Shanghai Jiao Tong University and the Key Laboratory of System Control and Information Processing, Ministry of Education of China, Shanghai, 200240, China, email: \{zhouyuanqiang,dwli,ygxi\}@sjtu.edu.cn.}
		\thanks{This work was supported by National Natural Science Foundation of China (NSFC) under Grant No. 61333009, 61521063, 61473317, 61573239. The work of Y. Zhou is supported by the China Scholarship Council.}}

	%\title{A Secure Control Learning  Framework for Cyber-Physical Systems under Sensor Attacks}
	%\author{Yuanqiang Zhou$^1$,~\IEEEmembership{Member,~IEEE,}
	%	Zhong-Ping Jiang$^1$,~\IEEEmembership{Fellow,~IEEE,}\\
	%	Kyriakos~G.~Vamvoudakis$^2$,~\IEEEmembership{Senior Member,~IEEE,}
	%	Wassim~M.~Haddad$^3$,~\IEEEmembership{Fellow,~IEEE}
	%	\thanks{$^1$Y.~Zhou and Z-P.~Jiang are with the Control and Networks Lab, Department of Electrical and Computer Engineering, Tandon School of
	%		Engineering, New York University, Brooklyn, NY 11201, USA, email: yz4734@nyu.edu, zjiang@nyu.edu.}
	%	\thanks{$^2$K.~G.~Vamvoudakis is with the Kevin T. Crofton Department of Aerospace and Ocean Engineering, Virginia Tech, Blacksburg, VA, 24061-0203, USA, e-mail: kyriakos@vt.edu.}
	%	\thanks{$^3$W.~M.~Haddad is with the School of Aerospace Engineering, Georgia Institute of Technology, Atlanta, GA
	%		30332-0150, USA, email:wm.haddad@aerospace.gatech.edu.}
	%	\thanks{This work was supported in part by NATO under grant No. SPS G5176, by ONR Minerva under grant No. N00014-18-1-2160, by an NSF CAREER,  by NAWCAD under grant No. N00421-16-2-0001 and by the Air Force Office of Scientific
	%		Research under Grant FA9550-16-1-0100.}}
	
	%%%%%%%%%%%%%%%%%%%%%%%%%%%%%%%%%%%%%%%%%%%%%%%%%%%%%%%%%%%%%%%%%%%%%%%%%%%%%%%
	
	\maketitle
	\thispagestyle{empty}
	\pagestyle{empty}
	
	%%%%%%%%%%%%%%%%%%%%%%%%%%%%%%%%%%%%%%%%%%%%%%%%%%%%%%%%%%%%%%%%%%%%%%%%%%%%%%%%
	\begin{abstract}
This paper investigates the data-driven predictive control problems for a class of continuous-time industrial processes with completely unknown dynamics.  The proposed approach employs the data-driven technique to get the system matrices online, using input-output measurements. Then, a model-free predictive control approach is designed to implement the receding-horizon optimization and realize the reference tracking. Feasibility of the proposed algorithm and stability of the closed-loop control systems are analyzed, respectively. Finally, a simulation example is provided to demonstrate the effectiveness of the proposed approach.
	\end{abstract}

\section{Introduction}
As a practically effective approach, model predictive control (MPC), or receding horizon control, has attracted notable attention in the field of industrial process control \cite{Karamanakos2014}.  To deal with optimal control problems, MPC can allow for industrial processes uncertainties and constraints much more straightforwardly than other methods \cite{Kothare1996,Mayne2000}. The core of all model-based predictive algorithms is to use ``open-loop optimal control" instead of ``closed-loop optimal control" within a moving horizon \cite{XI2013}. It brings a lot of robustness and reliability to allow the controller have the ability to recognize the control process. But, the acquisition of knowledge of a priori model seriously affects those performances of MPC. The practically measured data comes from the complicated processes and its utilization into MPC will greatly facilitate the design procedures, avoiding the need for initially accurate dynamic models \cite{Qin2012}.

Nowadays,  many research works focus on the data-driven predictive control, that is applicable on-line both to regulation and tracking control problems \cite{Xu2018,Yang2015,Shah2016,Wang2007,Hou2011,Lauri2010,Zhou2015,Aangenent2005,Ge2012,Ge2008,Kadali2003,zhou2019periodic}. More earlier works can refer to \cite{Skelton2018} and \cite{Aangenent2005}. A Markov data-based LQG control algorithm is suggested in \cite{Skelton2018} and data-based optimal control based on the system's Markov parameters is provided in \cite{Aangenent2005}. Both of them utilize the prior measurements to design the predictive control, along with its implementation on-line. Many results (see, \cite{Wang2007,Hou2011,Lauri2010}) still need some knowledge of the system information to fit the structure of the process model. Other works, like \cite{Li2013c}, focus on the polytopic uncertain systems with unmeasurable system states and a convex hull needs to be known. In \cite{Kadali2003,zhou2020synthesis}, data-driven subspace approach is introduced to design the predictive controller and in \cite{Gorges2017,Shah2016}, reinforcement learning approach is used to reduce the model-based dependence on predictive controller design procedures. We remark that all the aforementioned predictive control methods are designed based on the accurate model, as well as adequate uncertainty description of the linear or non-linear plant of the processes. 
%They make full use of the data-based control techniques to design the data-driven predictive controller. 

In this paper, the limitations of them are circumvented. We use the adequately measured data from the complicated industrial processes following the methods presented in \cite{Jiang2012,Jiang2017} and \cite{zhou2020secure,zhou2019secure}. To be precise, the data-driven learning technique of \cite{Jiang2012} will be employed to iteratively approximate the dynamical parameters, without requiring the \emph{a prior} knowledge of the system matrices. Then, the linear plant's version of \cite{Zhou2018} will be applied to predict the future trajectories, by following the continuous-time predictive control approach of \cite{Yang2015} but removing the assumption on partial knowledge of the system dynamics. Under this framework, the data-driven predictive control input can be generated on-line and can be used for the control of time-varying or nonlinear plants, since the algorithm is able to adapt to the actual dynamics by obtaining a linear model of the system at each sample. 
The contributions of the paper are three-fold. (i) A data-driven approach is proposed to adaptively approximate the system matrices, without requiring the \emph{a prior} knowledge of the system. (ii) A continuous-time data-driven MPC approach is developed for the continuous-time linear system, using repeatedly the state and input formation on some fixed time intervals. 
(iii) By implementing the proposed data-driven predictive control algorithm, both recursive feasibility of the optimization problem and closed-loop stability of the whole system are guaranteed.

The rest of this paper is organized as follows. In Section II, the problems are briefly formulated. In Section III, the data-driven predictive control approach with completely unknown dynamics algorithm is presented, and the feasibility and stability analysis are conducted. In Section IV, we apply the proposed approach to the optimal control problem of two continuous stirred tank reactor (CSTR). Conclusions are given in Section V.

\textbf{Notation}:  Through this note, $\mathbb{R}$ denotes the set of real numbers, $|\cdot |$ represents the Euclidean norm for a vector and the induced norm for a matrix. For a $a \in \mathbb{R}$, $\mathbb{R}_{\ge a}$ denotes the interval $[a , \infty)$ and $\mathbb{R}_{> a}$ denotes the interval $(a , \infty)$. For any vector $x$, $x^{\textrm{T}}$ denotes its transpose and ${\| x \|^2_P}$ represents ${{x^{\textrm{T}}}Px}$ for a real symmetrix and positive defined matrix $P$. For a real symmetric matrix $A \in \mathbb{R}^{n \times n}$, $A \succ 0$ or $A \succeq 0$ means that $A$ is positive definite or semi-positive definite, $\lambda_{M}(A)$ and $\lambda_{m}(A)$ denote the maximum and minimum eigenvalue of $A$, respectively. $\otimes$ indicates the Kronecker product operator, $\text{vec}(\cdot)$ denotes the vectorization operator and $\text{vec}^{-1}(\cdot)$ denotes the converse vectorization operator, i.e., $\text{vec}(A) = [a_1^{\textrm{T}}, \ldots, a_m^{\textrm{T}}]$ and $\text{vec}^{-1} [a_1^{\textrm{T}}, \ldots, a_m^{\textrm{T}}] = A$, where $a_i \in \mathbb{R}^n$ are the columns of a matrix $A \in \mathbb{R}^{n \times m}$. A continuous function $\alpha: \mathbb{R}_{\ge 0}\rightarrow \mathbb{R}_{\ge 0}$ is said to be a $\mathcal{K}$ function if it is strictly increasing, and $\alpha(s)>0$ for $s>0$ with $\alpha(0) = 0$. A continuous function $\alpha(\cdot)$ is said to be a $\mathcal{K}_{\infty}$ function if it is a $\mathcal{K}$ function, and $\alpha(s) \rightarrow \infty$ for $s \rightarrow 0$. For any piecewise continuous function $u(\cdot): \mathbb{R}_{\ge 0} \rightarrow \mathbb{R}^m$, $u^{[\textrm{r}]}$ denotes the $r$th order derivatives of $u(\cdot)$.

\section{Problem formulation and preliminaries}
Consider a continuous-time industrial process control system with the form
\begin{align}
{\dot x}(t) =&  A{x}(t) + B  u (t) , \quad x(t_0) = x_0, \label{eq:1}\\
y (t) = &  C x(t), \label{eq:2}
\end{align}
where $t \ge {t_0}$, ${x} \in \mathbb{R}^{n}$ is the measurable state; $ u\in \mathcal{U} \subset \mathbb{R}^{m}$ is the control input; $y \in \mathbb{R}^q$ is the controlled output fully available for feedback control design. Assume that $\mathcal{U} \subset \mathbb{R}^{m}$ is a nonempty compact convex set and contains the origin as its interior point. $A \in \mathbb{R}^{n \times n}$ and $B \in \mathbb{R}^{n \times m}$ are unknown system matrices with $(A,B)$ controllable, $(A, C)$ observable, satisfying $|A| \le A_M$, $|B| \le B_M$.   

Some standard assumptions are made on \eqref{eq:1} and \eqref{eq:2}. Similar assumptions can be found in \cite{Jiang2012,Zhou2018,Yang2015} for solving (cooperative) output tracking problems.

\noindent\textbf{Assumption 1.} \label{as:3}
The input constraint has the box constrain form as $\underline{u}\le u\le\overline{u}$ with elementwise inequality and $\underline{u},\overline{u}$ the respective lower and upper bounds.

\noindent\textbf{Assumption 2.} \label{as:1}
	There exists a constant matrix $K_0$ such that $A-B K_0$ is a Hurwitz matrix with $ - K_0 x(t) \in \mathcal{U}$.   

\noindent\textbf{Assumption 3.} \label{as:2}
	The input relative degree (IRD) of system \eqref{eq:1} is defined as $\rho$.

\noindent\textbf{Remark 1.} \label{re:1}
	In Assumption~\ref{as:3}, we refer the set $\mathcal{U}$ as a box constrain, which accurately describes nearly any set of standard mechanical actuators. Assumption~\ref{as:1} is made such that the initially feasibility can be achieved for the system \eqref{eq:1}. Assumption~\ref{as:2} imposes IRD for the system \eqref{eq:1}, which is used for simplify the solving of regulation problem \cite{Yang2015}. 

%Let ${\hat u}_{k}^*(t)$ denote the error between the optimal predictive controller at time $t = t_k$ and the one at the previous time $t = t_{k -1}$. 

For \eqref{eq:1} and \eqref{eq:2}, the output of the system should track the given reference $y_{\textrm{d}}(t), t\ge 0$. The tracking error can be given as $e(t) = y(t) - y_{\textrm{d}}(t), t\ge 0$. To that end, a sampled-data MPC, which is based on the repeated solution of an open-loop optimal control problem, is provided in this paper. At each time instant $t= t_k$, the state $y(t_k)$ is measured and then, the controller predicts the system behavior in the future over a prediction horizon $T$ by minimizing a certain objective cost function. The procedure is repeated at every sampling time instant $t_k$ for $k = 1, 2, \ldots$ 

For the system \eqref{eq:1} and \eqref{eq:2}, the cost function $J(t_k): = J(x(t_k), y_{\textrm{d}}(\cdot), \hat u_k(\cdot) )$, at time $t_k$, is defined as 
\begin{align}\label{eq:3-1}
J(x(t_k), y_{\textrm{d}}(s), \hat u_k(s) ) =& \int_{t_k }^{ t_k+ T} L(x(t_k), y_{\textrm{d}}(s), \hat u_k(s)) \textrm{d} s \notag \\
& + F(y_{\textrm{d}}(t_k + T), \hat y(t_k + T))
\end{align}
where, $L(\cdot, \cdot, \cdot)$ and $F(\cdot, \cdot)$ denote the stage and terminal cost functions with the from
\begin{align*}
L(x(t_k), y_{\textrm{d}}(t), \hat u_k( t))  =  \|e (t)\|_Q^2  +\| \hat u_k (t) ) \|_R^2, \quad t \ge 0, 
\end{align*}
where $Q  = Q^T \succ 0$ and $R = R^T \succeq 0$ are symmetric and sign definite weight matrices.

Then, at time $t_k$, the optimal control signal ${\hat u}_{k}^\star(s)$, $s\in [t_k, t_k + T]$, is obtained by solving the following finite-horizon optimal control problem as
\begin{subequations} \label{eq:4}
	\begin{align}
	&  \min\limits_{\hat u_k(s) \in \mathcal{U} } J(x(t_k), y_{\textrm{d}}(s), \hat u_k(s)  ) \notag \\
	&s.t. \quad {\dot {\hat x}}(s)  = \mathcal{H} \left(\hat x(s), \hat u_k(s) \right) {\Theta}   \label{eq:4A} \\
	& \qquad ~ \hat y (s) = C \hat x(s), {\hat x}(t_k) = x(t_k)   \label{eq:4B}   \\
	& \qquad~ {\hat u}(s) \in \mathcal{U}, \quad s\in [t_k, t_k +T ]   \label{eq:4C}
	\end{align}
\end{subequations}
where $\mathcal{H} (\cdot, \cdot): \mathbb{R}^n \times \mathbb{R}^m \rightarrow \mathbb{R}^{n \times ({n^2 + mn})} $ is defined as
\begin{align}\label{eq:5}
\mathcal{H}(x, u) = & \left[(x \otimes I_n)^{\textrm{T}} ~~  ( u \otimes I_n  )^{\textrm{T}} \right]
\end{align} 
and $\Theta$ denotes the vector of the dynamical parameters for (1), defined as
\begin{align}\label{eq:6}
\Theta =\left[ \text{vec}(A)^{\textrm{T}}  ~~ \text{vec}(B)^{\textrm{T}} \right]^{\textrm{T}}   \in \mathbb{R}^{n^2 + nm}
\end{align}
where $\text{vec}(\cdot)$ denotes the vectorization operator, defined as for any $ P= [p_{ij}] \in \mathbb{R}^{n \times l} , l \in \{n, m\}$, $\text{vec} (P) = \hat P$, where  
\begin{align}  \label{eq:7}
\hat P = [p_{11}, p_{21}, ..., p_{n1}, p_{12}, ..., p_{n2}, ..., p_{nl}]^{\textrm{T}}
\end{align}

%With the event-triggered mechanism, the solving of the MPC formulation (3) is implemented only when an event occurs. It can be seen that the predictive control input ${\hat u}_k^*(s)$ is an open-loop for $t \in [t_k, t_k + T_p]$,  since it only depends on the state $x(t_k)$ at each sampling time instant $t_k$. However, we can measure the error between the current output, $y(t)$, and the predicted state vector, $\hat y(t)$. By referring to this error, the triggering time instants are determined as
%\begin{equation}\label{eq8}
%t_{k+1} = \mathop {\inf }\limits_{\tau  \in [t_k,  t_k +T_p] }  \left\{  {\tau: {\left\| { \hat y(\tau) - y (\tau)} \right\|}_Q^2 \ge \sigma} \right\}
%\end{equation}
%where $\sigma$ is called as the triggering level to be designed. It is note that if (5) holds for some time $s \in [t_k, t_k + T_p]$, then take $t_{k + 1} = s $ ; if there do not exist time $ s \in [t_k. t_k + T_p] $ satisfying (5), then take $t_{k + 1} = t_k + T_p$.

%Since the derived optimization problem (3) is linear, by linear optimal control theory (Lewis Syrmos, 1995), when both $A, B$ and $c$ are accurately known, the solution to this problem can be found by solving the following well-known algebraic Riccati equation (ARE)
%\begin{align}\label{eq:4}
%- \dot P  = A^T P + PA + C^T Q C - P BR^{-1} B^T P
%\end{align}
%If (6) has a unique symmetric positive definite solution $P^*$, then the optimal predictive control law is given by
%\begin{align}\label{eq:5}
%u(t) = -K^* x(t) = R^{-1} B^T P^* x(t)
%\end{align}
%for $t \in [t_k, t_k + T_p]$. 

If the controller has the exact knowledge of system matrices $A$ and $B$, then we can directly solve the finite-time optimal control problem (4) by \cite{Frank2012} to get the desired predictive control policy ${\hat u}_k^\star$. But, in this paper, we seek a data-driven approach to remove the model-based dependence of (6) on the predictive controller design. Thus, the primary objective of this paper is to explore an online learning-based approach to find the data-driven predictive control policy, without requiring any knowledge of the system dynamics, just by using  \emph{a prior}  data of inputs and outputs. 

\section{Data-driven predictive control with completely unknown dynamics}
In this section, to facilitate the predictive controller design, the state ${\hat x}$ and the system parameters $\Theta$ of the predicted model (4a) are both estimated from the input-output measurements, using a simultaneous state and parameter estimator. Our proposed online-learning strategy does not rely on $\Theta$, i.e., either $A$ or $B$, which is totally data-driven approach. 

\subsection{Dynamical Parameters Estimator}

To obtain the dynamical parameters $\Theta$ in (6), i.e., $A$ and $B$, we present our online learning strategy with the previous input-output data. By using more prior data for continuous-time system (1)  than the instants of solving (4), we try to seek a least-squares optimization solution to get the estimations of all the dynamical parameters $\Theta$.  %First, we assume a stabilizing $K_0$ is known. 
To that end, at each time instant $t = t_k$, with enough prior date collected, we can find a minimal periodic time interval $\delta t$, such that for some $t_j< t_k$ of $j \in \{1, 2, ..., k-1\}$, there exists a constant integer $N_k$ satisfying $N_k \delta t = t_k - t_j$. 

To begin with, for any time instant $t \le  t_k$, doing the integrals at the time interval $[t - \delta t, t ]$ along the trajectory of (1) leads to
\begin{align}\label{eq:8}
x(t) - x(t - \delta t)  = & A \int_{t - \delta t}^{t}  x \textrm{d} \tau +  B \int_{t - \delta t}^{t }  \hat u^\star  \textrm{d} \tau  \notag \\
= & \mathcal{H} \left( \int_{t - \delta t}^{t}  x \textrm{d} \tau, \int_{t - \delta t}^{t }  \hat u^\star  \textrm{d} \tau   \right) {\Theta}  
\end{align}

By rearranging (8), we have the linear error system in the form of
\begin{align}\label{eq:9}
\mathcal{F}(t)  = \mathcal{G}(t) \Theta, ~~~ \forall t  \in \mathbb{R}_{[t_j, t_k]}
\end{align}
where the matrices $\mathcal{F} (\cdot): \mathbb{R}_{\ge 0} \rightarrow \mathbb{R}^n $ and $\mathcal{G} (\cdot): \mathbb{R}_{\ge 0} \rightarrow \mathbb{R}^{n \times (n^2 + m^2)} $ are defined,
\begin{align} 
\mathcal{F}(t) = & \begin{cases}   x(t) - x(t - \delta t), \quad t \in [\delta t, \infty)  \\
0, \qquad \qquad \qquad \quad t<\delta t
\end{cases}   \label{eq:10} \\
\mathcal{G}(t) = & \left[(\mathcal{A}(t)  \otimes I_n)^\textrm{T} ~~ (\mathcal{B}(t)  \otimes I_n)^\textrm{T} \right]    \label{eq:11}
\end{align}
where $\otimes$ denotes the Kronecker product,  and the vectors $ \mathcal{A} (\cdot) : \mathbb{R}_{\ge 0} \rightarrow \mathbb{R}^n$ and $\mathcal{B} (\cdot): \mathbb{R}_{\ge 0} \rightarrow \mathbb{R}^n $ are defined,
\begin{align}
\mathcal{A}(t) = & \begin{cases}   \int_{t - \delta t}^{t } x(\tau ) \textrm{d} \tau, \quad t \in [\delta t, \infty)  \\
0, \qquad \qquad \quad \quad t<\delta t
\end{cases}    \label{eq:12}  \\
\mathcal{B}(t) = & \begin{cases}   \int_{t - \delta t}^{t } \hat u^\star(\tau )  \textrm{d} \tau, \quad  t \in [\delta t, \infty)  \\
0, \qquad \qquad \quad \quad  t<\delta t
\end{cases}    \label{eq:13}
\end{align}

By (12)-(13), we note that for any time instant $t  \le t_k$, by using the previous input-output measurements $x(\tau)$ and $u(\tau)$, $\tau \le t$, $\mathcal{F}(t)$ and $\mathcal{G}(t)$ in (10)-(11) are both available. Then, by taking periodic sampling $\delta t$ for (1) at each time $t = t_k - i \delta t $  with $ i \in \{0, 1, 2, ..., N_k\}$, such that $ t_j \le t \le t_k $, then we can get a set of composite data, 
\begin{align}\label{eq:14}
\mathfrak{D}_k =   \bigcup \{\mathcal{F}_{i}, ~  \mathcal{G}_{i} \}_{i =0}^{N_k}
\end{align}
which satisfies 
\begin{align}\label{eq:15}
\mathcal{F}_i = \mathcal{G}_{i} \hat \Theta + \epsilon , ~ \forall i \in \{1, 2, ..., N_k\} 
\end{align}
where $\mathcal{F}_i = \mathcal{F} (t_k -i\delta t)$, $\mathcal{G}_i = \mathcal{G} (t_k -i\delta t)$, $\hat \Theta  = [ \hat A^T  ~~ \hat B^T ]^T$ is an estimate of $ \Theta$ in (6) and $\epsilon$ is the estimation error due to the data-driven approximation of $ \Theta$ by using $\mathfrak{D}_k$.

To bring $\epsilon$ to its minimum value, we have the following optimization problem (OP): 
\begin{align}  \label{eq:16}
&  \min\limits_{\hat \Theta } ~~ \epsilon^\textrm{T} \epsilon  \\
&s.t. \quad (15) ~ \text{and} ~ (9)-(13)  \notag
\end{align}

Then, an adaptive least-square method is used here to give a solution to problem (16). But, before that, the full rank of dataset $\mathfrak{D}_k$ is defined as follows. To guarantee the existence of the solution of (16), we also give the full rank condition for dataset $\mathfrak{D}_k$.

\noindent\textbf{Definition 1. }\label{de:1}
	At each time instant $t_k$ for some integer $N_k$, the data stack $\mathfrak{D}_k$ is said to have full rank, if there exists an integer $N_0 >0$, such that for all $N_k \ge N_0$, we have the matrix $\mathfrak{A}_k$, defined as 
	%		\begin{align*}
	%\mathfrak{A}_k = \sum_{i = 1}^{N_k} \mathcal{G}_{i} ^T \mathcal{G}_{i}  \in \mathbb{R}^{(n^2 + nm) \times (n^2 + nm) }
	%	\end{align*} 
	$\mathfrak{A}_k : = \sum_{i = 1}^{N_k} \mathcal{G}_{i} ^\textrm{T} \mathcal{G}_{i}  \in \mathbb{R}^{(n^2 + nm) \times (n^2 + nm) }$, satisfying
	\begin{align}\label{eq:17}
	\text{rank}(\mathfrak{A}_k ) = n^2 + nm 
	\end{align}

Then, we have the following lemma  to give a sufficient condition to guarantee the full rank property of dataset $\mathfrak{D}_k$. The proof is similar to \cite{Ioannou1996}, and here it is omitted. 

\noindent\textbf{Lemma 1.}\label{le:1}
	A data stack $\mathfrak{D}_k$ has full rank, if there exists a constant $\underline{d} > 0$, such that \begin{align}\label{eq:18}
	0 < \underline{d} < \gamma_{m}(\mathfrak{A}_k)
	\end{align} 

Assume that the measurements of the inputs and the outputs are prior collected at an enough large number $N_k \gg n^2 + nm$ points of time $t_k - i \delta t$ with $i \in \{0, 1, ..., N_k\}$, which makes data stack $\mathfrak{D}_k$ has full rank by Lemma~\ref{le:1}. Then, at time $t = t_k$, one can use the data stack $\mathfrak{D}_k$ to evaluate the unknown dynamical parameters $\Theta$ in (6). By solving (16), the following learning-based update law is obtained
\begin{align}\label{eq:19}
\dot{\hat{\Theta}} =  \eta_{\theta} \sum\limits_{i =1}^{N_k} \mathcal{G}_i^\textrm{T} \left(\mathcal{F}_i - \mathcal{G}_i \hat{\Theta}  \right)
\end{align} 
where $0 < \eta_{\theta} \in \mathbb{R}$ is a constant learning rate.

For ease of exposition, the original continuous-time system (1) is expressed in the form of,
\begin{align}\label{eq:20}
{\dot x}(t)  = \mathcal{H} \left(x(t), u(t) \right) {\Theta}
\end{align} 
By considering the linear error system with parameters update law (19), it follows that 
\begin{align}\label{eq:21}
{\dot x}(t)  = & \mathcal{H} \left( x(t) , u(t)  \right)  \hat  {\Theta} + w(x(t), u(t))   \notag \\
: = & \tilde {A} x(t) + \tilde B  u(t)  + w(t)
%\dot{\hat{\Theta}} = & ~ \eta_{\theta} \sum\limits_{i =1}^{N_k} \mathcal{G}_i^T \left(\mathcal{F}_i - \mathcal{G}_i \hat{\Theta}  \right)
\end{align} 
where %$\hat \Theta =[ \hat A^T  ~~ \hat B^T ]^T$, $ \tilde A = \text{vec}^{-1}(\hat A)$, 
$ \tilde B = \text{vec}^{-1}(\hat B)$,  $ \tilde A = \text{vec}^{-1}(\hat A)$, $ \text{vec}^{-1}(\cdot)$ denotes the converse vectorization operator, that is, for any vector $\hat P \in \mathbb{R}^{n\times l} $ in (7), we have $  \text{vec}^{-1}(\hat P)= \tilde P = [p_{ij}] \in \mathbb{R}^{n\times l} $. And $w(t):  = w(x(t), u(t))$ is the continuous approximation error resulting from parameters uncertainty (19). If $\epsilon = 0$, then $\hat \Theta =\Theta$, it implies $\tilde A = A$ and $\tilde B = B$, thus we have $w(x, u) = 0$. 

The following theorem analyses the property of the term $w$ with respect to (21). 

\noindent\textbf{Theorem 1.}\label{th:1}
	The approximate error $w(t)$ in (20) is slowly time-varying, bounded and satisfies $\lim_{t \rightarrow \infty}  w(t) =0 $. \\
\emph{Proof. }
	By considering the closed-loop dynamics (20) and (21), we can refer to $w(t)$ as the unknown disturbance caused by the parameters uncertainty of (19). Then, $w(t)$ satisfies
	\begin{align}\label{eq:22}
	w(t) = \mathcal{H} \left(x(t), u(t) \right) \left( \hat{\Theta} -  \Theta  \right)
	\end{align} 
	
	Letting $ \tilde{\Theta}  = \hat{\Theta} -  \Theta $ and bringing (19) to (22) lead to
	\begin{align}\label{eq:23}
	|\dot w(t)| \le  &| \mathcal{H} \left(x(t), u(t) \right) |  |\dot{ \tilde{\Theta}} | \notag \\
	\le & \eta_{\theta}   | \mathcal{H} \left(x(t), u(t) \right) |  | \sum\limits_{i =1}^{N_k} \mathcal{G}_i^\textrm{T}  \mathcal{G}_i |   | \dot{ \tilde{\Theta}} |  \notag \\
	= &  \eta_{\theta}   | \mathcal{H} \left(x(t), u(t) \right) |  | \mathfrak{A}_k|   | \dot{ \tilde{\Theta}} |  
	\end{align} 
	
	By Lemma 1, we have $\gamma_{m} (\mathfrak{A}_k) \le  | \mathfrak{A}_k| \le \gamma_{M} (\mathfrak{A}_k)$. Then, by (19), it implies $\dot {\tilde{\Theta} } = - \eta_{\theta} \mathfrak{A}_k \tilde{\Theta}$ and thus $\lim_{t \rightarrow \infty} \tilde{\Theta} (t) =0 $ and $\tilde{\Theta}  \le  e^{ - \eta_{\theta}  \gamma_{m} (\mathfrak{A}_k) }  \Theta(t_k - N_k \delta t)$. So, $ | \dot{ \tilde{\Theta}} | \le \eta_{\theta}  \gamma_{M} (\mathfrak{A}_k) e^{ - \eta_{\theta}  \gamma_{m} (\mathfrak{A}_k) }  \Theta(t_k - N_k \delta t): = C_{\Theta}$. Besides, for a fixed dataset $\mathfrak{D}_k$ in (14), we have $ | \mathcal{H} \left(x(t), u(t) \right)|$ in (23) bounded for some real constant $M \in \mathbb{R}$ at each time $t = t_k -i \delta t$. Thus, we have
	\begin{align}\label{eq:24}
	|\dot w(t)| \le  \eta_{\theta}   M C_{\Theta}  \gamma_{M} (\mathfrak{A}_k)  
	\end{align} 
	It implies that $w(t)$ is slowly time-varying and $\lim_{t \rightarrow \infty}  w(t) =0 $. By (23),(24), we also have $|w(t)| \le e^{ \eta_{\theta}   M C_{\Theta}  \gamma_{M} (\mathfrak{A}_k) } $, which means that $w(t)$ is bounded.

\noindent\textbf{Remark 2. }
	Note that in \eqref{eq:21}, the conventional receding horizon expression $A x(t) + Bx(t)$ depending on the unknown matrices $A, B$ is replaced by the term $ \mathcal{H} \left(\hat x(s), \hat u_k(s) \right) \tilde{\Theta}$, where $\tilde{\Theta}$ can be obtained by repeatedly learning from the states and inputs measurements. Furthermore, this learned results will not affects the convergence of the system by Lemma~\ref{th:1}. Therefore, \eqref{eq:19}  plays an important role in identifying the system dynamics from the \emph{a prior} data. As a result, the requirement of the system matrices in predicting the behavior of \eqref{as:2} can be replaced by the state and input information measured online.

\subsection{Receding-horizon Optimization}
To facilitate the data-driven predictive controller design for the system (21) with the dynamical parameters estimator (19), the receding-horizon predictive control problem of (4), at time instant $t = t_k$, can be reformulated,
\begin{subequations}  \label{eq:25}
	\begin{align}
	& {\hat u}_k^\star(s) = \arg \mathop {\min }\limits_{{{\hat u}}(t) \in \mathcal{U}} {J}({{x}}({t_k}),y_\textrm{d}(s), {{\hat u_k}}(s)) \notag \\
	&s.t. \quad {\dot {\hat x}}(s) = \tilde A {\hat x} (s) + \tilde B \hat u_k (s),  \label{eq:25A} \\
	& \qquad ~ \hat y (s) = C \hat x(s),  {\hat x}(t_k) = x(t_k),  \label{eq:25C}  \\
	& \qquad~ {\hat u}(s) \in \mathcal{U}, \quad s \in [t_k, t_k+ T].  \label{eq:25E}
	\end{align}
\end{subequations}

To solve the optimization problem (25),  under Assumption~\ref{as:2}, we define the decision variables as ${\bar u}_k (s)= [\hat u_k^\textrm{T}(s), (\hat u_k^{[1]})^\textrm{T}(s), \ldots, (\hat u_k^{[r]})^\textrm{T}(s)]$ for some  control order $r$ larger than $\rho\ge 1$. Note that the first term of ${\bar u}_k (\tau)$ is the to-be-optimized control input $\hat u_k (\tau)$ in (25). More generally, for the control law $\hat u_k (\tau) $ with a large enough control order $r$, we let $\hat u_k^{[l]} (\tau) =0 $ for any integer $l \ge r$. 

Then, for the output prediction of optimization problem (25), by following (2), the future output $y(s) = y(t + \tau), t = t_k, k = 1, 2, \ldots$, in the moving horizon $\tau \in [0, T]$ is approximated by Taylor series expansion,
\begin{align}\label{eq:26}
y(t + \tau) = y(t) + \tau y^{[1]}(t) + \cdots + \frac{\tau^{r}}{ r!} y ^{[r] } (t) + O(\tau^{r})
\end{align}
where the $i$-th derivative of the output $y^{[i]}(t)$ with $i \in \{1, 2, ..., \rho, ..., r\}$ is obtained by
\begin{align}
{y^{[i]}} = & C{{\tilde A}^i}x + \sum\limits_{k = 0}^{i - 1} C {{\tilde A}^{i - 1 - k}}{w^{[k]}},i = 1,...,\rho  - 1  \label{eq:27}\\
{y^{[j]}} = & C{{\tilde A}^j}x + \sum\limits_{k = 0}^{j - \rho } C {{\tilde A}^{j - 1 - k}}\tilde B{u^{[k]}}  \notag \\
& + \sum\limits_{k = 0}^{j - 1} C {{\tilde A}^{j - 1 - k}}{w^{[k]}}, j = \rho, ..., r \label{eq:28}
\end{align}

By rewriting the output $y(t + \tau)$ in a compact form, it follows that
\begin{align}\label{eq:29}
y(t + \tau )  = \begin{bmatrix}  T_1 (\tau) & T_2 (\tau) \end{bmatrix}   \begin{bmatrix}
Y_1 \\  Y_2 \end{bmatrix}
\end{align}
where ${{T}}_1 (\tau) = \left[ 1, \tau, \ldots,  \frac{\tau^{\rho -1}}{(\rho - 1)!}  \right] $, $T_2(\tau) = \left[    \frac{\tau^{\rho}}{ \rho !}, \ldots,  \frac{\tau^{ r }}{ r!}  \right]$, ${Y}_1=\left[ y^\textrm{T}, (y ^{[1]})^\textrm{T} , \ldots,  (y ^{[\rho - 1 ]}) ^\textrm{T}  \right]^\textrm{T}$, ${Y}_2=\left[  (y^{[\rho]} ) ^\textrm{T}, (y^ {[\rho + 1 ]})^\textrm{T} , \ldots,  (y ^{[r]})^\textrm{T} \right]^\textrm{T}$, and 
\begin{align}
Y_1 = & \mathcal{A}_1 x +  \mathcal{B}_1 \bar w   \label{eq:30} \\
Y_2 =&  \mathcal{A}_2 x +  \mathcal{B}_2 \bar w  +  \mathcal{B}_3 \bar u \label{eq:31}
\end{align}
where ${\bar w} (s)= \left[ w^\textrm{T}(s),  (w^{[1]})^\textrm{T}(s), \ldots,   (w^{[r]})^\textrm{T} (s) \right]^\textrm{T}$ and $\mathcal{A}_1, \mathcal{B}_1, \mathcal{A}_2, \mathcal{B}_2$ and $ \mathcal{B}_3$ are defined,
\begin{align}\label{eq:32}
\mathcal{A}_1 =& \begin{bmatrix}  C \\ C \tilde A   \\ \vdots \\  C \tilde A^{\rho - 1}  \end{bmatrix} ,    \mathcal{B}_1 =\begin{bmatrix}  C &  0 & \cdots & 0  \\ C \tilde A &  C & \cdots & 0  \\ \vdots & \vdots &  \ddots &  \vdots \\  C \tilde A^{ \rho -1 } & C \tilde A^{ \rho -2}  \tilde B & \cdots & C \end{bmatrix}  \notag \\
\mathcal{A}_2 = &  \begin{bmatrix}  C \tilde A^{\rho } \\ C\tilde A^{\rho +1 }    \\ \vdots \\  C \tilde A^{ r }  \end{bmatrix}, 
\mathcal{B}_2 \!=\! \begin{bmatrix}  C\tilde A^{\rho-1 }  & 0 \cdots & 0 \\ C\tilde A^{\rho }  &  C \tilde A^{\rho-1 }  & \cdots & 0  \\ \vdots & \vdots & \ddots &  \cdots \\  C \tilde A^{r  }  & C \tilde A^{ r -1}  & \cdots & C\tilde  A^{ \rho -1 }  \end{bmatrix}  \notag  \\
\mathcal{B}_3= & \begin{bmatrix}  C \tilde A^{\rho-1 } \tilde B & 0 &  \cdots & 0 \\ C\tilde A^{\rho }  \tilde B&  C \tilde A^{\rho-1 } \tilde B  & \cdots & 0  \\ \vdots & \vdots & \ddots & \cdots \\  C {\tilde A}^r  \tilde B & C {\tilde  A}^{r -1}  \tilde B & \cdots & C{ \tilde A}^{ \rho -1 } \tilde B  \end{bmatrix}
\end{align}

Then, for the reference signal $y_d(t + \tau)$, by Theorem 1, we have $w^{[k]} (\tau)\approx 0$ for $k=1, 2, \ldots, r$. Thus, $y_d(t + \tau)$ satisfies (29)-(31) with $Y_{1,d} = \left[ y_d^\textrm{T}, (y_d ^{[1]})^\textrm{T} , \ldots ,  (y_d ^{[\rho - 1 ]}) ^\textrm{T}  \right]$ and $ Y_{2, d} =\left[  (y_d^{[\rho]} ) ^\textrm{T}, (y_d^ {[\rho + 1 ]})^\textrm{T} , \ldots ,  (y _d^{[r]})^\textrm{T} \right]$. 
But, along with the solution of optimization problem (25) without considering the parameters uncertainty from (21), it leads to
\begin{align}
& \hat y(t + \tau )  = \begin{bmatrix}  T_1 (\tau) & T_2 (\tau) \end{bmatrix}   \begin{bmatrix}
Y_1 \\  Y_2    \end{bmatrix}  \label{eq:33} \\
& Y_1 = \mathcal{A}_1 x, \quad Y_2 =  \mathcal{A}_2 x + \mathcal{B}_3 {\bar u}  \label{eq:34}\\
& \hat u(t + \tau) = T (\tau) \bar u  \label{eq:35}
\end{align}
where $T(\tau) = [T_1(\tau)~ T_2(\tau) ]$.

Then, by bring (34)-(35) into the given performance index (3), one can go to,
\begin{align}\label{eq:36}
J(t_k)  = & \int_{0}^{T}   \left[ \tilde Y_1^\textrm{T}, \tilde Y_2^\textrm{T} \right] \begin{bmatrix}
\Xi_1^\textrm{T} (\tau) \\  \Xi_2^\textrm{T}(\tau)  \end{bmatrix}  \left[ \Xi_1 (\tau) ,  \Xi_2 (\tau) \right]  \begin{bmatrix}
\tilde Y_1 \\ \tilde Y_2 \end{bmatrix}    \notag  \\
& + \int_{0}^{T}  \hat { u} ^\textrm{T}  T_3^\textrm{T} (\tau) R  T_3 (\tau)   \bar { u}   \textrm{d}\tau + F(\tilde Y_i (t_k + T))
\end{align}
where $\tilde Y_i =Y_i -Y_{i,d}$ for $i \in \{1, 2\}$, and $\Xi_i (\tau)= \sqrt{Q}  T_i(\tau) $ for $i \in \{1, 2\}$. 
By defining ${\mathcal{T}}_{i, j} = \int_{0}^{T} \Xi_i^\textrm{T} \Xi_j \textrm{d}\tau$ with $i, j \in \{1, 2\}$ and ${\mathcal{T}} = \int_{0}^{T}  T^\textrm{T} R T \textrm{d}\tau$, we note that ${\mathcal{T}}_{1, 2} = {\mathcal{T}}_{2, 1}$. Thus, the performance index (37) can also be rewritten as 
\begin{align}\label{eq:37}
J(t_k)  = & \tilde Y_1^\textrm{T} {\mathcal{T}}_{1,1} \tilde Y_1 + 2 \tilde{Y}_1^\textrm{T} {\mathcal{T}}_{1,2} \tilde {Y}_2  +  \tilde{Y}_2^\textrm{T} {\mathcal{T}}_{2,2} \tilde{Y}_2 \notag \\
& +   {\bar u} ^\textrm{T}  \mathcal{T}   {\bar u}   +  F(\tilde Y_i (t_k + T))
\end{align}

Due to equation (34),  taking partial derivative of $J(t_k)$ with respect to $ {\bar u}_k$  yields 
\begin{align}\label{eq:38}
\frac{{\partial J}}{{\partial \bar u}} =& 2  \left(  (\frac{{\partial \tilde Y_2}}{{\partial \bar u}})^\textrm{T} {\mathcal{T}}_{2, 2}  \tilde Y_2  + {\mathcal{T}}_4 \right)  {\bar u}  + 2 (\frac{{\partial \tilde Y_2}} {{\partial \bar u}}  )^\textrm{T}  {\mathcal{T}}_{1,2}^\textrm{T}  \tilde{Y}_1     \notag  \\
= & 2  \left(   {\mathcal{B}}_3^\textrm{T}  {\mathcal{T}}_{2,2} {\mathcal{B}}_3 +  {\mathcal{T}}_4  \right)  {\bar u}  + 2 {\mathcal{B}}_3^\textrm{T}  {\mathcal{T}}_{1,2}^\textrm{T} \tilde{Y}_1  \notag   \\
& + 2   {\mathcal{B}}_3^\textrm{T}  {\mathcal{T}}_{2,2} \left( \mathcal{A}_2 x - Y_{2,d}\right)
\end{align}

By letting  ${\partial J}/{{\partial {\bar u}_k  }} = 0$, we can get the optimized predictive control law ${\bar u}_k^\star$ as
\begin{align}\label{eq:39}
{\bar u}_k^\star=& -  \left(   {\mathcal{B}}_3^\textrm{T}  {\mathcal{T}}_{2,2} {\mathcal{B}}_3 +  {\mathcal{T}}_4  \right) ^{-1}   {\mathcal{B}}_3^\textrm{T}  \left( {\mathcal{T}}_{2,2} \bar Y_2 - {\mathcal{T}}_{1,2}^\textrm{T} \tilde{Y}_1 \right)
\end{align}
where $\bar Y_2 =  Y_{2,d} -  \mathcal{A}_2 x $. Then, taking the first row of the optimized control law (39), the continuous-time predictive control law, applied to the plant, is given by
\begin{align}\label{eq:40}
\hat u_k^\star (t) = I_u   {\bar u}_k^\star
\end{align}
where $I_u = [1, 0, ..., 0]_{1 \times (r+1)}$. 

\noindent\textbf{Remark 3.}
	Note that from (40), the existence of the optimized solution ${\bar u}_k^\star$ depends on the reversibility of matrix $\mathcal{M} = {\mathcal{B}}_3^\textrm{T}  {\mathcal{T}}_{2,2} {\mathcal{B}}_3 +  {\mathcal{T}}_4$, with ${\mathcal{B}}_3$ computed from (32), $\tilde A$ and $\tilde B$  calculated from (19).  Thus, before we implement the receding-horizon optimization, we first check the reversibility of matrix $\mathcal{M}$ by removing the repeated columns of dataset $\mathfrak{D}_{k}$, only left the distinct columns. 

\subsection{Handing constraints}

To deal with optimal control problems, MPC can allow for industrial processes uncertainties and constraints much more straightforwardly than other methods {[}2{]}, {[}3{]}. Assumption~\ref{as:3} will allow us to use a specialized active-set algorithm which is more efficient and easier to implement. A box-constraint solver can be immediately generalized to any linear inequality constraints using slack variables. In the following, we formalize two classical ways to enforce the control limits. 

\subsubsection{Saturating Functions}

A conventional attempt to enforce box constraints is to clamp the controls in the forward-pass. The element-wise clamping, or projection operator, is denoted by $sat\text{(\ensuremath{\cdot})}$,which is the input saturation function defined as
\begin{align}
sat(u) & =\left[sat(u_{1})\:sat(u_{2})\:\cdots\:sat(u_{m})\right]\label{eq:3}\\
sat(u_{i}(t)) & =\begin{cases}
u_{i}(t) & {\textstyle if}\,u_{i,min}<u_{i}(t)<u_{i,max}\\
u_{i,min} & {\displaystyle if}\,u_{i}(t)\le u_{i,min}\\
u_{i,max} & if\,u_{i}(t)\ge u_{i,max}
\end{cases}\nonumber 
\end{align}
with $u=[u_{1}\:u_{2}\cdots u_{m}]$, and $u_{i,min}\le0$ and $u_{i,max}\ge0$
are the boundaries of $i$th control input of system (1).

It implies that it is tempting to simply replace the obtained control in the forward-pass with 
\[
\hat{u}=sat(u^{*})
\]
However, the corresponding search direction may not be a descent direction anymore, harming convergence. 

\subsubsection{Squashing Functions}

Another way to enforce box constraints is to introduce a sigmoidal squashing function $s(u)$ on the controls
\[
x_{i+1}=f(x_{i},s(u_{i}))
\]
where $s(\cdot)$ is an element-wise sigmoid with the vector limits
\[
\underset{u\to-\infty}{lim}s(u)=\underline{u},\quad\underset{u\to\infty}{lim}s(u)=\overline{u}
\]
For example, $s(u)=\frac{\overline{u}-\underline{u}}{2}tanh(u)+\frac{\overline{u}+\underline{u}}{2}$
is such a function. A cost term should be kept on the original $u$ and not only on the squashed $s(u)$. Otherwise it will reach very
high or low values and get stuck on the plateau. An intuition for the poor practical performance of squashing is given by the nonlinearity of the sigmoid. Since the backward pass uses a locally quadratic approximation of the dynamics, significant higher order terms will always have a detrimental effect on convergence.

\subsection{Our Proposed Algorithm}
Our proposed method is summarized as the following Algorithm $1$.

\smallskip
\begin{algorithm}
	\caption{Data-driven predictive control algorithm}\label{al:1}
	\KwData{current period $t$; a initially stable control $u^\star = -K_0 x$; prediction horizon $T$; terminal cost $F$;}
	\KwResult{Optimal MPC input $u^{\ast}_{t}$}
	$t_k = k \leftarrow 0$\;
	Collect data $\mathfrak{D}_k$ in (14) \;
	\While{}{Collect the data and form $\mathfrak{D}_k$\;
		\If{(18) is satisfied}{Generate the estimator (19) by using $\mathfrak{D}_k$\;
			Implement (25) to obtain the optimized control (40)\;
			Time evolves continuously with $t$\;}
		$k\leftarrow k+1$\;}
\end{algorithm}

\subsection{Performance Analysis}

Before proceeding further, we first introduce the following definition and lemma.

\noindent\textbf{Definition 2.}  For the system (21), given a compact set $\mathbb{E}$, with $\{0\}   \subset \mathbb{E}\subseteq \mathbb{R}^n$ and $\mathbb{E}$ being a robustly positively invariant set, if there exists a positive definite function $V(\cdot): \mathbb{R}^n \to \mathbb{R}_{\ge 0}$, such that,
	\begin{align}
	V( x) &\ge \alpha_1 (| x|), ~~V( x) \le \alpha_2 (| x|)  + c_1   \label{eq:41}\\
	V(\dot { x}) & -V( x)  \le - \alpha_3 (| x|) + \alpha_ 4(|w|) +c_2  \label{eq:42}
	\end{align}
	for all $t \in \mathbb{R}_{\ge 0} $ with $\alpha_1, \alpha_2, \alpha_3$ being $\mathcal{K}_{\infty}$ function, $\alpha_4$ being $\mathcal{K}$ function, and $c_1, c_2 \ge 0$. Then, the function $V(\cdot)$ is a regional input-to-state practical stability (ISpS)-type Lyapunov function in $\mathbb{E}$ for the system.  

Based on Definition 2, we can have the following lemma, directly borrowed from \cite{Yang2015}. 

\noindent\textbf{Lemma 2.}  Given a robust positively invariant set $\mathbb{E}$ for the system (21), if it admits an ISpS-type Lyapunov function $V(\cdot)$, then the system is regional ISpS in $\mathbb{E}$, and all the signals of the closed-loop system with the control input $u_{\mathbb{E}}$ are bounded, where  $u_{\mathbb{E}}(t)$ denotes the control such that the set $\mathbb{E}$ is an invariant region satisfying the constraints.

Note that by Assumption 1, the developed optimization problem of data-driven predictive control in (25) is initially feasible with $u_0 = -K_0 x$, then the global stability can be proved by using the Lemma 2.

\noindent\textbf{Theorem 2.}\label{th:2}
	Suppose the Assumption 1 and 2 hold for system with a robust positively invariant set $\mathbb{E}$, then the closed-loop control  plant (1) under the continuous-time MPC law (40), is globally asymptotically stable.

\emph{Proof.}
	The proof is composed by two parts, feasibility and convergence.
	
	\emph{Feasibility:} Consider any time $t_k$ such that the problem of (25) has a solution and the optimal input $\hat u_{k}^\star$ is implemented for time $[t_k, t_{k+1})$.  Assumed that at $t_{k+1}$, $\hat y(t_{k+1}) = y(t_{k+1})$. 	Therefore, the remaining piece of optimal input $\hat u_k^\star (s), \tau \in [t_{k+1}, t_{k} + T]$ satisfies the input constraints. Thus, we construct the control input as,
	\begin{align}\label{eq:43}
	\hat u_{k+1}(\tau) = \begin{cases}  \hat u_{k}^\star(\tau),  ~~ \tau \in [t_{k+1}, t_{k} + T]  
	\\u_{\mathbb{E}}(\tau), ~~ \tau \in [t_{k} + T, t_{k+1} + T]  
	\end{cases}
	\end{align}
	where $u_{\mathbb{E}}(\tau)$ makes the desired reference reached and the constraints satisfied. Thus, the predictive control problem is feasible at $t_{k+1}$. It implies that feasibility of the problem at $t_k$ implies the recursive feasibility at $t_{k+1}$. 
	
	\emph{Convergence:} Let the optimal cost function at $t_k$ as the value function $V(x(t_k)) = J^\star(x(t_k), y_\textrm{d}(s), \hat u_k^\star(s) )$. If $V(x(t_k))$ is strictly decreasing, the tracking error $e$ will converge to the origin. To this end, we write the value function at $t_k$ as,
	\begin{align}\label{eq:44}
	V(t_k) =& \int_{t_k}^{t_k+ T} \left( \|e(\tau)\|_Q^2 +\| \hat u_{k}^\star(\tau) ) \|_R^2\right) \textrm{d} \tau  \notag \\
	& + F (y_d(t_k + T), \hat y(t_k + T)) 
	\end{align}
	
	Then, by applying $\hat u_{k+1}(t)$  in (44) to the system, beginning from $y(t_{k+1})$, one has
	\begin{align}\label{eq:45}
	J(t_{k+1}) =& \int_{t_{k+1}}^{t_{k+1}+ T} \left( \|e(\tau)\|_Q^2 +\| \hat u_{k+1}(\tau) ) \|_R^2\right) \textrm{d} \tau  \notag \\
	& + F (y_\textrm{d}(t_{k+1} + T), \hat y(t_{k+1} + T)) 
	\end{align}
	
	By substituting (45) in (46), it follows,
	\begin{align}\label{eq:46}
	J(t_{k+1}) =& V(t_k) - \int_{t_k}^{t_{k+1}} \left( \| e(\tau)\|_Q^2 +\| \hat u_k^\star(\tau) ) \|_R^2\right) \textrm{d} \tau  \notag \\
	& -   F (y_\textrm{d}(t_k + T), \hat y(t_k + T))       \notag \\
	& +   \int_{t_k + T}^{t_{k+1}+ T} \left( \|e (\tau)\|_Q^2 +\| \hat u_{k+1}(\tau) ) \|_R^2 \right) \textrm{d} \tau   \notag \\
	& + F (y_\textrm{d} (t_{k+1} + T), \hat y(t_{k+1} + T)) 
	\end{align}
	
	Note that $\mathbb{E}$ is a robust positively invariant set, it implies that for all $x \in \mathbb{E}$, we have
	\begin{align}\label{eq:47}
	(\partial F / \partial x) (Ax + Bu_{\mathbb{E}}) + L (x, y_\textrm{d}, u_{\mathbb{E}})  \le 0
	\end{align}
	where $L (x, y_\textrm{d}, u_{\mathbb{E}}) = \|e \|_Q^2 +\| u_{\mathbb{E}}  \|_R^2 $. By integrating the inequality (48) along the trajectory of $\dot {\hat x} = \tilde A \hat x + \tilde B u_{\mathbb{E}}$ and $\hat y = C\hat x$, we have $\int_{t_k + T}^{t_{k+1}+ T} \left( \|e (\tau)\|_Q^2 +\| u_{\mathbb{E}}(\tau) ) \|_R^2 \right) \textrm{d} \tau \le F (y_\textrm{d}(t_k + T), \hat y(t_k + T)) - F (y_\textrm{d} (t_{k+1} + T), \hat y(t_{k+1} + T)) $. Thus, by (47), one gets,
	\begin{align}\label{eq:48}
	V(t_k) - J(t_{k+1}) \! \le \! - \int_{t_k}^{t_{k+1}} \left( \| e(\tau)\|_Q^2 +\| \hat u_k^\star (\tau) ) \|_R^2\right) \textrm{d} \tau
	\end{align}
	
	Further, based on $V(t_{k+1}) = J^\star (t_{k+1}) $, we obtain, 
	\begin{align}\label{eq:49}
	V(t_k) - V(t_{k+1}) \! \le \! - \int_{t_k}^{t_{k+1}} \left( \| e(\tau)\|_Q^2 +\| \hat u_k^\star(\tau) ) \|_R^2\right) \textrm{d} \tau
	\end{align}
	It implies that $V(t_{k+1})$ is strictly decreasing. Hence, the proof is complete.

\newcounter{MYtempeqncnt}
\begin{figure*}[!t]
	% ensure that we have normalsize text
	\normalsize
	% Store the current equation number.
	\setcounter{MYtempeqncnt}{\value{equation}}
	% Set the equation number to one less than the one
	% desired for the first equation here.
	% The value here will have to changed if equations
	% are added or removed prior to the place these
	% equations are referenced in the main text.
	\setcounter{equation}{49}
	\begin{align}
	\label{eqn_dbl_x}
	A =  \begin{bmatrix} 
	-17.98 & -295.866 & 0    &   0 &   0 &    0\\
	0.0207  &  0.1889 &   0.0704  &       0      &   0      &   0\\
	0  &  0.3879  &  0.8000       &  0     &    0     &    0\\
	0.0977   & 0     &    0       &      -18.01 & -295.87 &  0\\
	0  &  0.0617   & 0      &        0.0131   & 0.0433  &  0.0589\\
	0  &  0   &      0      &        0    &     0.3787&   -0.622 \end{bmatrix}, \quad 
	B   =  \begin{bmatrix} 
	17.8996   &   -13.781  \\
	-0.0131   &  0.0101  \\
	0   &    0  \\
	17.8636& 17.8636  \\
	0.0082 &  0.0082 \\
	0 0 \end{bmatrix}
	%	C \!\!  = \!\! \begin{bmatrix} 
	%	0 & 0\\ 362.995 & 0\\ 0& 0 \\0 &
	%	0\\ 0& 362.995 \\0 & 0 \end{bmatrix}^T 
	\end{align}
	% Restore the current equation number.
	\setcounter{equation}{\value{MYtempeqncnt}}
	% The IEEE uses as a separator
	\hrulefill
	% The spacer can be tweaked to stop underfull vboxes.
	\vspace*{4pt}
\end{figure*}

% The previous equation was number five.
% Account for the double column equations here.
\addtocounter{equation}{1}

\section{Application to two-CSTR process}\label{sec5}
Consider two continuous stirred tank reactor (CSTR) system with a full description in \cite{Cao2002,Wang2016}. The open-loop model is a six-state continuous model. The system matrices $A$ and $B$ are directly taken from \cite{Cao2002}, described in the form of (1), as you can see in (\ref{eqn_dbl_x}). The system output variables $y_1 = 362.995 x_2$ and  $y_2 = 362.995 x_4$, denoting the two tank outlet temperatures. The control problem is to maintain the two tank temperatures at desired values $y_d(t) =[y_{1d}(t)~y_{2d}(t)]^\textrm{T} $, where $y_{1d}(t) = 10$ when $0 \le t < 5$s and $y_{1d}(t) = 7$ when $ t \ge 5$s, $y_{2d}(t) = 10$ when $0 \le t < 5$s and $y_{2d}(t) = 4$ when $ t \ge 5$s. The constraints is,
\begin{align*}
\mathcal{U} = \{u =[u_1~ u_2]^\textrm{T}: |u_1| \le 80 , ~~|u_2| \le 70 \}
\end{align*}

In order to illustrate the efficiency of the proposed approach, the precise knowledge of $A$ and $B$ is not used in the design of the predictive controllers. Since the physical system is not stable, the initial stabilizing feedback gain is set as $K_0$,
\begin{align*}
K_0 = \begin{bmatrix} 
\begin{smallmatrix}
-4.8949 & -3426.8 & -158.1712    &   -0.0320 &   -43.7963 &   -1.4675 \\
0.1  & 0  &   86.2934  &       1.1730      &   2.3886      &   104.8756
\end{smallmatrix}
\end{bmatrix}
\end{align*}

The weighting matrices $Q$ and $R$ are set to be $Q = diag([10~ 100~ 10~ 10~ 100~ 10])$ and $R = diag([1~ 1])$, respectively. In the simulation, the initial values are selected at the origin. The state and input information is collected over each interval of $0.01$s. When time arrives at $t =2$ s, all the inputs and outputs are repeatedly used to approximate the matrices $A$ and $B$ with $\eta_\theta = 0.85$. The predictive control also starts at $t = 2$ s with the prediction horizon $T = 1$ s. Since then, the control input is immediately updated by solving the problem (25), and the convergence of of $A_k: = \tilde{A}$ and $B_k: = \tilde{B}$ to their actual values is attained after 10 iterations. The procedure of solving (25) is repeated over a fixed interval of $0.1$s. 
The convergence of $A_k$ and $B_k$ to their actual values is illustrated in Fig.1. The trajectories of the output variables and the flow rates are shown in Fig.2. 
It can be seen that the data-driven predictive control algorithm can stabilize the system, without requiring the system matrices. 

\begin{figure}
	\centering
	\includegraphics[width=3.6in]{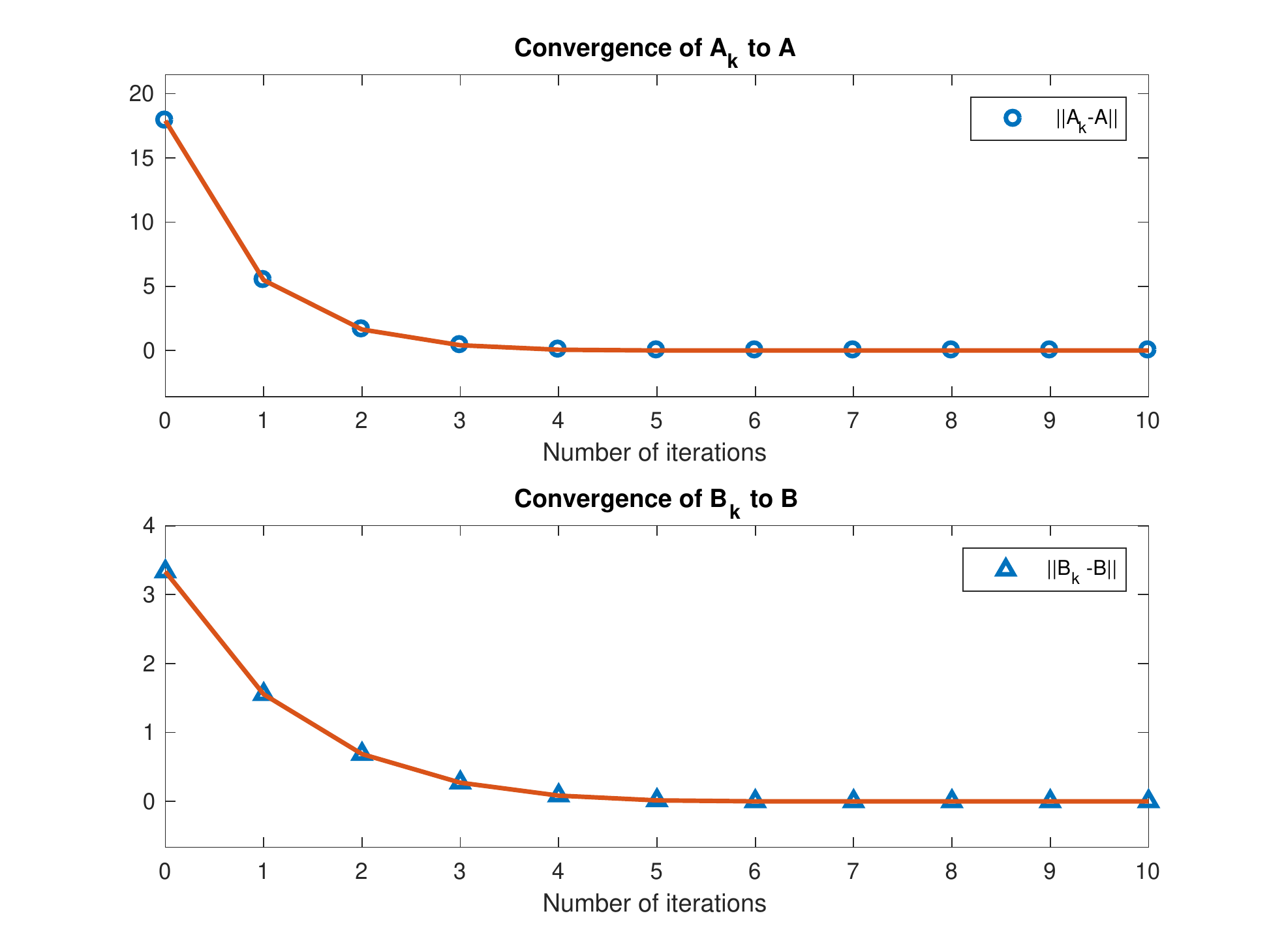}\\
	\caption{Convergence of $A_k$ and $B_k$ to their actual values during the control process.}\label{fig:3}
\end{figure}
\begin{figure}
	\centering
	\includegraphics[width=3.6in]{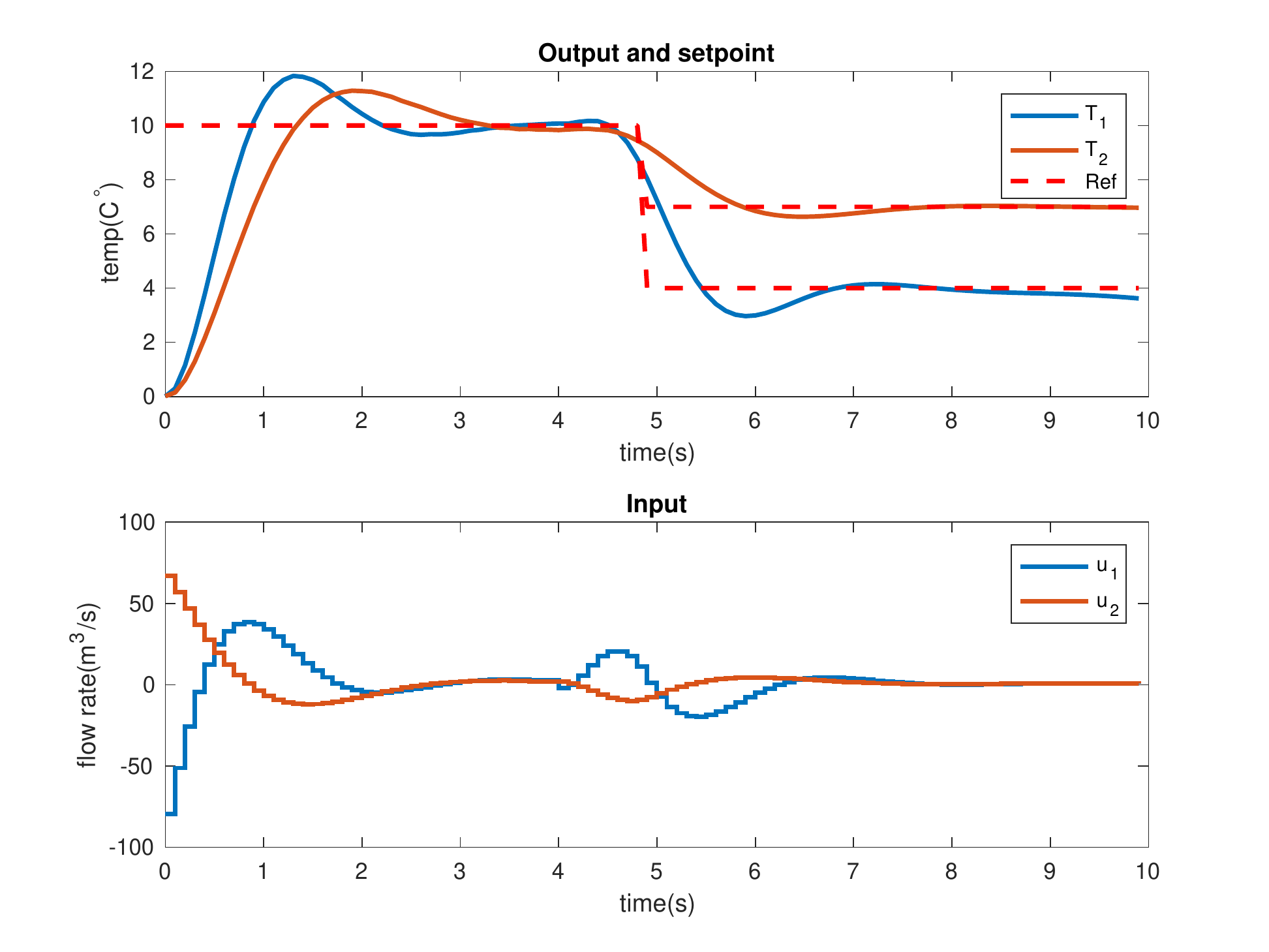}\\
	\caption{The trajectories of the output variables and the flow rates.}\label{fig:2}
\end{figure}

\section{Conclusion}\label{sec10}
In this paper, a data-driven predictive control approach for continuous-time linear system with completely unknown dynamics has been provided. This method solves the infinite-horizon optimal control problem, using the system inputs and outputs information collected online, without knowing the system matrices. The methodology developed in this paper may serve as a computational tool to study the finite-horizon adaptive optimal control of uncertain nonlinear systems. Some related work has appeared in \cite{Ge2008}, which was developed using neural networks, and also in our recent work \cite{Zhou2018}, which proposed a framework of distributed MPC to handle the asynchronous communication, using the \emph{a prior} information  associated with the interconnected neighbors to a distributed optimal design.

%	\hrule
%	\vspace{0.25cm}
%	%{\bf Algorithm 1:} 
%	
%	\vspace{0.25cm}
%	\hrule
%	\begin{algorithmic}[1]
%		%\Procedure{}{}
%		\Require A initially stable control $u^\star = -K_0 x$; prediction horizon $T$; terminal cost $F$;
%		\State Set $t = t_k = k = 0$;
%		\While { Time arrives at $ t = t_k$}
%		\Repeat
%		\State Collect data $\mathfrak{D}_k$ in (14);
%		\Until {(18) is satisfied}
%		\State	Generate the estimator (19) by using $\mathfrak{D}_k$.
%		\State Implement (25) to get the optimized control (40);
%		\State Apply the control (40) to the system (1);
%		\State Time evolves continuously with $t$; 
%		\EndWhile
%		\vspace{0.25cm}
%		\hrule
%	\end{algorithmic}
%\end{algorithm}

%\begin{ack}
%Place acknowledgments here.
%\end{ack}
%\bibliographystyle{Bibliography/IEEEtranTIE}
\bibliographystyle{IEEEtran}
%\bibliography{Bibliography/IEEEabrv,CDC}\ %IEEEabrv instead of IEEEfull
\bibliography{RefLib}\ %IEEEabrv instead of IEEEfull

%\section{A summary of Latin grammar}    % Each appendix must have a short title.
%\section{Some Latin vocabulary}              % Sections and subsections are supported  
                                                                         % in the appendices.
\end{document}